\input amstex.tex 
\documentstyle{amsppt} 
\NoRunningHeads 
\magnification=\magstep1 
\baselineskip=12pt 
\parskip=5pt 
\parindent=18pt 
\topskip=10pt 
\leftskip=0pt 
\pagewidth{32pc} 
\pageheight{45pc} 
\topmatter 
\title 
Automorphisms of Finite Order 
on Gorenstein del Pezzo Surfaces 
\endtitle 
\author D. -Q. Zhang 
\endauthor 
\address 
\newline \noindent 
Department of Mathematics, National University of Singapore, 
2 Science Drive 2, 
\newline \noindent 
Singapore 117543, Republic of Singapore
\endaddress 
\par \noindent 
\email 
\newline \noindent 
matzdq$\@$math.nus.edu.sg 
\endemail 
\dedicatory 
\enddedicatory 
\subjclass 
Primary 14J50; Secondary 14J26 
\endsubjclass 

\abstract 
In this note we shall determine all actions of groups of prime 
order $p$ with $p \ge 5$ on Gorenstein del Pezzo (singular) surfaces $Y$
of Picard number 1. We show that every order-$p$ element in
$\text{\rm Aut}(Y)$ ($= \text{\rm Aut}({\widetilde Y})$,
${\widetilde Y}$ being the minimal resolution of $Y$)
is lifted from a projective transformation of ${\bold P}^2$.
We also determine when $\text{\rm Aut}(Y)$ is finite
in terms of $K_Y^2$, $\text{\rm Sing} Y$ and the number 
of singular members in $|-K_Y|$.
In particular, we show that either $|\text{\rm Aut}(Y)| = 2^a3^b$ for some
$1 \le a+b \le 7$, or for every prime $p \ge 5$ there is at least
one element $g_p$ of order $p$ in $\text{\rm Aut}(Y)$ (hence $|\text{\rm Aut}(Y)|$ is infinite).
\endabstract 
\endtopmatter 

\document 
\head Introduction 
\endhead 

We work over the complex numbers field ${\bold C}$. 
In this note we study the pair $(X, G)$ of a normal rational surface $X$ 
and a finite group $G$ of automorphisms on $X$. 

\par 
When $X$ is smooth, this subject had been studied by S. Kantor more than 
one hundred years ago [K]. 
It was continued by Segre, Manin, Iskovskih, Gizatullin
and many others [S], [M1, 2], [I], [G]. 
See also [H1], [H2].
In [DO], the group of automorphisms of any general del Pezzo surface is described
and it turns out that its discrete part is equal to the
kernel of the Cremona representation on the moduli space of $n$ points in ${\bold P}^2$.
Very recently, de Fernex [dF] constructed all the Cremona transformations of ${\bold P}^2$ 
of prime order, where he employed the methods different from those used 
by Dolgachev and the author of this note in [ZD].

\par
In [ZD], minimal pairs $(X, G)$ with prime order $p = |G|$ was considered. 
In particular, using the recent Mori theory, 
it was shown there that if the $G$-invariant sublattice 
of $\text{\rm Pic} X$ has rank 1 then $p \le 5$ unless $X = {\bold P}^1$; 
the short and precise classification of these pairs, modulo equivariant
isomorphism, was also given there. 
In [MM] and [MZ3], more general situation was considered where $X$ may be 
singular or even open. 

\par 
In this note we will consider the case where $Y$ is a Gorenstein del Pezzo 
singular surface of Picard number 1. So the pair $(Y, G)$ is 
automatically minimal in the generalized sense of [ZD]. 
In contrast to the smooth case in [ZD], we shall show that 
every prime number $p \ge 5$ is realizable as the order of some
element in $\text{\rm Aut}(Y)$ for some $Y$.
This actually confirms the common belief that for a family
of algebraic surfaces the automorphism group of a special
member should be bigger than that of a generic member;
see Remark D below.
See also [O] for the case of family of K3 surfaces.

\par
The other contrasting result is that in smooth case
there are minimal pairs $(X, {\bold Z}/(5))$ with $X \ne {\bold P}^2$ so
the action of ${\bold Z}/(5)$ on $X$ is not lifted from
a regular action on ${\bold P}^2$ [ZD, Table 1], while in
singular case, every element of prime order $p \ge 5$ in $\text{\rm Aut}(Y)$ 
($= \text{\rm Aut}({\widetilde Y})$) is lifted from a
projective transformation of ${\bold P}^2$ (Theorem A).

\par
Moreover, we show that for a given $Y$, 
the inclusion ${\bold Z}/(p_0) \subseteq \text{\rm Aut}(Y)$
for one single prime $p_0 \ge 5$ will guarantee the 
inclusion ${\bold Z}/(p) \subseteq \text{\rm Aut}(Y)$
for every prime $p \ge 5$ (Theorem C).

\par 
We begin with a definition. 
A normal projective surface $Y$ is a {\it Gorenstein del Pezzo surface} 
if $Y$ has only rational double singularities 
and if the anti-canonical divisor 
$-K_Y$ is ample. 

\par
As we see from the classification of higher dimensional
algebraic varieties (see [KM]), a minimal (resp. canonical) model will have some 
terminal (resp. canonical) singularities. In surface case, canonical singularities
are just rational double singularities (or Du Val singularities,
or Dynkin type ADE singularities, or rational Gorenstein singularities in other notation).
Also Gorenstein del Pezzo surfaces $Y$ appear naturally as degenerate fibres in a family of
(smooth) del Pezzo surfaces (= Fano varieties of dimension 2);
the minimal model conjecture (still unsolved for dimension 4 or greater)
claims that every algebraic variety is birational to either a minimal
model or a variety with a Fano fibration (whose singular fibres are
varieties with some mild singularities and ample anti-canonical divisor).
This is also the motivation why we study singular del Pezzo surfaces.

\par
Gorenstein del Pezzo singular surface $Y$ of Picard number 1 satisfies 
$1 \le K_Y^2 \le 9$ and $K_Y^2 \ne 7$ (see e.g. [MZ1]). 
The case $K_Y^2 = 8$ (resp. $K_Y^2 = 9$) occurs if and only if 
$Y$ is the quadric cone ${\overline \Sigma}_2$ in ${\bold P}^3$ (resp.
$Y = {\bold P}^2$).
We now state our main results.

\par \vskip 1pc \noindent 
{\bf Theorem A.} {\it Let $p \ge 5$ be a prime number. 
Let $Y$ ($\ne \overline{\Sigma}_2$)
be a Gorenstein del Pezzo (singular) surface of Picard number $1$.
Then we have:}

\par \noindent
(1) {\it Modulo equivariant isomorphism, there is either none or only one
or exactly $p+1$ of non-trivial ${\bold Z}/(p)$-action(s) on $Y$.}

\par \noindent
(2) {\it Each non-trivial ${\bold Z}/(p)$-action  on $Y$
equals, modulo equivariant isomorphism, one of those in Examples $1.1-1.9$.}

\par \noindent
(3) {\it Every order-$p$ element $g$ in $\text{\rm Aut}(Y)$
($= \text{\rm Aut}({\widetilde Y})$, ${\widetilde Y}$
being the minimal resolution of $Y$) is lifted from a projective
transformation $\overline{g}$ of ${\bold P}^2$,
i.e., there is a birational morphism $\mu : {\widetilde Y} \rightarrow {\bold P}^2$
such that $\mu(g y) = \overline{g} \mu(y)$ for all $y$ in ${\widetilde Y}$.}

\par \vskip 1pc \noindent 
{\bf Theorem B.} {\it Let $p \ge 5$ be a prime number. 
Let $Y$ be a Gorenstein del Pezzo (singular) surface of Picard number $1$.} 

\par \noindent 
(1) {\it If $K_Y^2 = 6$, then there are exactly $p+1$ of
non-trivial ${\bold Z}/(p)$-actions 
on $Y$, modulo equivariant isomorphism.}

\par \noindent 
(2) {\it If $K_Y^2 = 5$, then there is a unique 
non-trivial ${\bold Z}/(p)$-action 
on $Y$, modulo equivariant isomorphism.} 

\par \noindent 
(3) {\it If $K_Y^2$ is in $\{3, 4\}$, 
then there is either only one or exactly $p+1$
of non-trivial ${\bold Z}/(p)$-action(s) 
on $Y$, modulo equivariant isomorphism.} 

\par \noindent 
(4) {\it Suppose $K_Y^2 = 2$. Then there is a unique 
(resp. there is no any) non-trivial
${\bold Z}/(p)$-action 
on $Y$, modulo equivariant isomorphism,
if $\text{\rm Sing} Y \ne A_7$ (resp. $\text{\rm Sing} Y = A_7$).} 

\par \noindent 
(5) {\it Suppose $K_Y^2 = 1$. Then there is a unique (resp. there is no any)
non-trivial ${\bold Z}/(p)$-action
on $Y$, modulo equivariant isomorphism, if $|-K_Y|$ has exactly two
(resp. more than two) singular members (see {\bf 2.5}).} 

\par \vskip 1pc \noindent 
{\bf Theorem C.} {\it Let $Y$ be a Gorenstein
del Pezzo (singular) surface of Picard number $1$. }

\par \noindent
(1) The following are equivalent:

\par \noindent
(1a) {\it $\text{\rm Aut}(Y)$ has finite order.}

\par \noindent
(1b) {\it $|\text{\rm Aut}(Y)| = 2^a 3^b$ with $1 \le a+b \le 7$.}

\par \noindent
(1c) {\it Either $\text{\rm Sing} Y = A_7$ (and hence $K_Y^2 = 2$), or
$K_Y^2 = 1$ and $|-K_Y|$ contains at least three singular members
(i.e., $K_Y^2 = 1$ and $Y$ is different from those four in Examples $1.1 - 1.4$).}

\par \noindent
(2) {\it The following are equivalent:}

\par \noindent
(2a) {\it $\text{\rm Aut}(Y)$ has infinite order.}

\par \noindent
(2b) {\it $\text{\rm Aut}(Y)$ contains an element of prime order $p_0$ with $p_0 \ge 5$.}

\par \noindent
(2c) {\it For each prime number $p \ge 5$, there is at least one element $g_p$
of order $p$ in $\text{\rm Aut}(Y)$.}

\par \vskip 1pc \noindent
{\bf Remark D.} We like to compare Theorems B and C with known results for smooth $Y$. 

\par \noindent
(1) It is known that the order of $\text{\rm Aut}(Y)$ of a generic rational surface
(the blow up of ${\bold P}^2$ at very generic points)
with $K_Y^2 \le 5$ is a factor of $5!$ (see [DO] or [Ko, Main Theorem]).  

\par \noindent
(2) Let $Y$ be a (smooth) del Pezzo surface with $K_Y^2 = 3$ or $4$.
It is known that there is at most one (resp. no any) ${\bold Z}/(p)$-action
on $Y$, modulo equivariant isomorphism, 
if $p = 5$ (resp. $p \ge 7$ is prime); see [S, pp. 147-152],
[H1, Theorem 1.1] and [H2, Theorem 5.3].

\par
An interesting observation made by I. Dolgachev is that
the quotient surface of the degree 6 Del Pezzo surface modulo
an involution, is a 4-nodal Segre cubic surface.

\par \noindent
(3) The following can be deduced from [DO] or the main Theorem in [Ko]. 
Let $Y$ be a generic rational surface.
If $K_Y^2 \le 4$, then there is no prime order $p$ ($p \ge 3$) 
automorphism of $Y$. Suppose $K_Y^2 = 5$; then there is a unique
(resp. no any) non-trivial ${\bold Z}/(p)$-action 
on $Y$, modulo equivariant isomorphism, if $p = 5$ (resp. if $p \ge 7$ is prime).
If $K_Y^2 = 6, 7$, then for each prime $p \ge 5$, 
there are exactly $p+1$ of ${\bold Z}/(p)$-actions 
on $Y$, modulo equivariant isomorphism. 

\par \noindent
(4) Let $Y$ be a Gorenstein del Pezzo singular surface of Picard number 1. 
[Y1] or [Y2] implies:
If $K_Y^2 \ge 2$, the $\text{\rm Sing} Y$ determines uniquely the isomorphism 
class of $Y$. For $K_Y^2 = 1$, see [Y1, 2] or Proposition 2.5.

\par \noindent 
(5) We believe that a similar classification is achievable for $p = 2, 3$ 
though the list will be much longer and less elegant. 

\par \noindent
(6) The case $K_Y^2 = 1$ has been considered in [Z2].

\head
Acknowledgement
\endhead
The author would like to thank Professor I. Dolgachev for the discussion
on generic del Pezzo surfaces and Professor M. Reid on maps
between singular del Pezzo surfaces.
This work was partially supported by an Academic Research Fund of 
National University of Singapore.

\head 
Section 1. Examples 
\endhead 

We begin with a definition.
A rational elliptic (smooth) surface $f : X \rightarrow {\bold P}^1$ is called 
{\it extremal} if  $f$ is relatively minimal (i.e., $K_X^2 = 0$) and 
if the Mordell Weil group $\text{\rm MW}(f)$ of all sections is torsion (cf. [MP], [OS]); 
the latter is equivalent to saying that the Picard lattice over {\bf Q} of $X$ 
is generated by a single section and fibre components. 

\par
We will first give Examples 1.1 - 1.5, where $p \ge 5$ is a prime number and 
$\zeta_p = exp(2 \pi \sqrt{-1}/p)$. 
We will define an action of $\langle g \rangle \cong {\bold Z}/(p)$ 
on certain extremal rational elliptic (smooth) surface $X$. 
In the Weierstrass equations for $X$ below 
(cf. [MP, Tables 5.1-5.3]), we use $v$ to denote the parameter of the base curve. 
In Examples 1.1-1.4, a different choice of section to be blown down by the 
map $X \rightarrow Y$ will result in isomorphic $Y$ because 
$\text{\rm MW}(f)$ acts (on $X$ and) transitively on the set of all sections. 
Also except Example 1.5, the action of $\langle g \rangle$ on $X$ 
stabilizes all negative curves (i.e., $(-2)$-curves 
in fibres and $(-1)$-curves as sections); see Theorems 2.1 and 2.2.

\par \vskip 1pc \noindent 
{\bf Example 1.1.} Let $f : X \rightarrow B \cong {\bold P}^1$ be the 
{\it unique} rational elliptic (smooth) surface with singular fibres of type $II$ and $II^*$. 
Its Weierstrass equation is given here: 
$$y^2z = x^3 + v^5z^3.$$ 
We define an order-$p$ automorphism $g$ on $X$ as follows: 
$$g^*: (x, y, z; v) \rightarrow (\zeta_p^{-5} x , y, \zeta_p^{-15} z; \zeta_p^6 v).$$ 
Then $g$ preserves the fibration $f$ and stabilizes exactly two fibres 
of type $II^*$ at $v = 0$ and of type $II$ at $v = \infty$. 
The restriction $g|B$ has order $p$. 

\par \vskip 1pc 
Clearly, $g$ stabilizes the unique section $E$ of $f$, which is given by 
$[x, y, z; v] = [0, 1, 0; v]$. We label the type $II^*$ 
fibre as $\sum_{i=1}^6 i C_i + 4C_7 + 2C_8 + 3C_9$ 
so that $E + \sum_{i=1}^8 C_i$ is a linear chain and $C_9 . C_6 = 1$. 
Let $X \rightarrow Y_k$ ($0 \le k \le 4$) be the composite of the smooth blow-down 
of $E+\sum_{i=1}^k C_i$ and the contraction of 
the $(-2)$-curves $\sum_{i=k+2}^9 C_i$. 
Then $Y = Y_k$ satisfies : 

\par \vskip 1pc \noindent 
(*) $Y$ is Gorenstein del Pezzo, 
Picard number $\rho(Y) = 1$, $K_Y^2 = k+1$ and \text{\rm Sing} $Y$ is one of 
$E_8, E_7, E_6, D_5, A_4$ (depending on $k$). 

\par \vskip 1pc \noindent 
The $g$ on $X$ induces a regular action of order-$p$ on $Y_k$, 
which we denote by the same letter $g$. 

\par \vskip 1pc \noindent 
{\bf Example 1.2.} Let $f : X \rightarrow B \cong {\bold P}^1$ be 
the {\it unique} rational elliptic (smooth) surface with singular fibres of type $III$ and $III^*$. 
Its Weierstrass equation is given here: 
$$y^2z = x^3 + v^3xz^2.$$ 
We define an order-$p$ automorphism $g$ on $X$ as follows: 
$$g^*: (x, y, z; v) \rightarrow (\zeta_p^{-3} x , y, \zeta_p^{-9}z; \zeta_p^4 v).$$ 
Then $g$ preserves the fibration $f$ and stabilizes exactly two fibres 
of type $III^*$ at $v = 0$ and of type $III$ at $v = \infty$. 
The restriction $g|B$ has order $p$. 
\par \vskip 1pc 
Clearly, $g$ stabilizes the only two sections $E, E'$ of $f$, which are given by 
$[x, y, z; v] = [0, 1, 0; v]$ and $[x, y, z; v] = [0, 0, 1; v]$. 
We label the type $III^*$ fibre as $\sum_{i=1}^4 iC_i + \sum_{j=5}^8 a_j C_j$ 
so that $E + \sum_{i=1}^7 C_i$ is a linear chain. 
Let $X \rightarrow Y_k$ ($k = 0, 1, 2$) be the composite of the smooth blow-down 
of $E + \sum_{i=1}^k C_i$ and the contraction of the $(-2)$-curves 
$\sum_{i=k+2}^8 C_i$ and the $(-2)$-curve in the type $III$ fibre 
not meeting $E$. 
Then $Y = Y_k$ satisfies : 

\par \vskip 1pc \noindent 
(*)  $Y$ is Gorenstein del Pezzo, 
Picard number $\rho(Y) = 1$, $K_Y^2 = k+1$ and \text{\rm Sing} $Y$ is one of 
$E_7 + A_1, D_6 + A_1$ and $A_5 + A_1$ (depending on $k$). 

\par \vskip 1pc \noindent 
The $g$ on $X$ induces a regular action of order-$p$ on $Y_k$, 
which we denote by the same letter $g$. 

\par \vskip 1pc \noindent 
{\bf Example 1.3.} Let $f : X \rightarrow B \cong {\bold P}^1$ be 
the {\it unique} rational elliptic (smooth) surface with singular fibres of type $IV$ and $IV^*$. 
Its Weierstrass equation is given here: 
$$y^2z = x^3 + v^4z^3.$$ 
We define an order-$p$ automorphism $g$ on $X$ as follows: 
$$g^*: (x, y, z; v) \rightarrow (\zeta_p^{-2} x , y, \zeta_p^{-6}z; \zeta_p^3 v).$$ 
Then $g$ preserves the fibration $f$ and stabilizes exactly two fibres 
of type $IV^*$ at $v = 0$ and of type $IV$ at $v = \infty$. 
The restriction $g|B$ has order $p$. 

\par \vskip 1pc 
Clearly, $g$ stabilizes the only three sections $E, E', E''$ of $f$, 
which are given by 
$[x, y, z; v] = [0, 1, 0; v]$ and $[x, y, z; v] = [0, \pm v^2, 1; v]$. 
We label the type $IV^*$ fibre as $\sum_{i=1}^3 iC_i + \sum_{j=4}^7 a_j C_j$ 
so that $E + \sum_{i=1}^5 C_i$ is a linear chain. 
Let $X \rightarrow Y_k$ ($k = 0, 1$) be the composite of the smooth blow-down 
of $E + \sum_{i=1}^k C_i$ and the contraction of the $(-2)$-curves 
$\sum_{i=k+2}^7 C_i$ and the two $(-2)$-curves in the type $IV$ fibre 
not meeting $E$. 
Then $Y = Y_k$ satisfies : 

\par \vskip 1pc \noindent 
(*)  $Y$ is Gorenstein del Pezzo, 
Picard number $\rho(Y) = 1$, $K_Y^2 = k+1$ and \text{\rm Sing} $Y$ is one of 
$E_6 + A_2$ and $A_5 + A_2$ (depending on $k$). 

\par \vskip 1pc \noindent 
The $g$ on $X$ induces a regular action of order-$p$ on $Y_k$, 
which we denote by the same letter $g$. 

\par \vskip 1pc \noindent 
{\bf Example 1.4.} Let $f : X_J \rightarrow B \cong {\bold P}^1$ be the rational elliptic 
(smooth) surface with two singular fibres of type $I_0^*$ and $J$-invariant of a general fibre 
equal to the constant $J = 4r^3/(4r^3 + 27 s^2)$. Here $r, s$ are 
in ${\bold C}$ so that $4r^3 + 27s^2 \ne 0$. 
Its Weierstrass equation is given as follows, 
where $v$ is the parameter of the base curve 
$$y^2z = x^3 + rv^2xz^2 + sv^3z^3.$$ 
We define an order-$p$ automorphism $g$ on $X$ as follows: 
$$g^*: (x, y, z; v) \rightarrow (\zeta_p^{-1} x , y, \zeta_p^{-3}z; \zeta_p^2 v).$$ 
Then $g$ preserves the fibration $f$ and stabilizes exactly the two fibres 
of type $I_0^*$ at $v = 0$ and $v = \infty$. 
The restriction $g|B$ has order $p$. 

\par \vskip 1pc 
Clearly, $g$ stabilizes the only four sections $E, E', E'', E'''$ of $f$, 
which are given by 
$[x, y, z; v] = [0, 1, 0; v]$ and $[x, y, z; v] = [x_i v, 0, 1; v]$, 
where $x_i$ are the roots of $x^3+rx + s = 0$. 
We label a type $I_0^*$ fibre as $\sum_{i=1}^2 iC_i + \sum_{j=3}^5 a_j C_j$ 
so that $E + \sum_{i=1}^3 C_i$ is a linear chain. 
Let $X_J \rightarrow Y_{J,k}$ ($k = 0, 1$) be the composite of the smooth blow-down 
of $E + \sum_{i=1}^k C_i$ and the contraction of the $(-2)$-curves 
$\sum_{i=k+2}^5 C_i$ and the four $(-2)$-curve in the another type $I_0^*$ fibre 
not meeting $E$. 
Then $Y = Y_{J, k}$ satisfies : 

\par \vskip 1pc \noindent 
(*)  $Y$ is Gorenstein del Pezzo, 
Picard number $\rho(Y) = 1$, $K_Y^2 = k+1$ and \text{\rm Sing} $Y$ is one of 
$2D_4$ and $D_4 + 3A_1$ (depending on $k$). 

\par \vskip 1pc \noindent 
The $g$ on $X_{J}$ induces a regular action of order-$p$ on $Y_{J,k}$, 
which we denote by the same letter $g$. 
\par \vskip 1pc \noindent 
In general, if $f : X \rightarrow B$ ($\cong {\bold P}^1$) has 
singular fibre type $I_0^*, I_0^*$, then $X = X_J$ for some $J$ [MP, Theorem 5.4]. 
We remark also that the isomorphism class of 
$Y_{J, k}$ depends (resp. does not depend) on $J$ when $k = 0$ 
(resp. $k = 1$); see Proposition 2.5. 

\par \vskip 1pc \noindent 
{\bf Example 1.5.} Let $f : X \rightarrow B \cong {\bold P}^1$ be 
the unique rational elliptic (smooth) surface with singular fibres of type $I_1, I_1, I_5, I_5$. 
Its Weierstrass equation is given here: 
$$y^2z = x^3 + Axz^2 + Bz^3,$$ 
where $A = -3(1 - 12v+14v^2+12v^3+v^4)$, $B = 2(1 - 18v + 75v^2 + 75v^4 + 18 v^5 + v^6)$. 
The Mordell Weil group $G = \text{\rm MW}(f) \cong {\bold Z}/(5)$ acts on $X$ 
naturally as translations (of general fibres). In particular, 
it acts transitively on the set of all 5 sections. 
The $G$ stabilizes every fibre. So the restriction $G|B$ is trivial. 

\par \vskip 1pc
In Examples 1.6 - 1.10 below, we let $Y$ be a Gorenstein del Pezzo 
singular surface of Picard number 1. Let ${\widetilde Y} \rightarrow Y$ be the minimal resolution 
and $D$ the exceptional divisor. 
Denote by $\#D$ the number of irreducible components of $D$. 
Then the Picard number $\rho({\widetilde Y}) = 1 + \#D$. Hence 
$K_Y^2 = K_{\widetilde Y}^2 = 9 - \#D$. For $K_Y^2 \ge 2$, the $\text{\rm Sing} Y$
determines uniquely the isomorphism class of $Y$ (see Proposition 2.5);
we will also use Figure 5 in [Y1] or [Y2, Ch 4] containing 
all negative curves on ${\widetilde Y}$, which are either $(-1)$-curves
or $(-2)$-curves.

\par \vskip 1pc \noindent 
{\bf Example 1.6.} Let $Y$ be the {\it unique} Gorenstein del Pezzo 
surface of Picard number 1 and with $\text{\rm Sing} Y = A_2+A_1$. Then $K_Y^2 = 6$. 

\par 
We first construct such unique surface $Y$. 
Let $P_1 = [0, 0, 1]$ and $P_2 = [0, 1, 0]$. 
Let $\mu: {\widetilde Y} \rightarrow {\bold P}^2$ 
be the blow up of $P_1$ and its 2 infinitely near 
points lying on the proper transform of the line $L_{P_1P_2}$, such that 
$\mu^{-1}(L_{P_1P_2}) = D_1 + E + D_3 + D_2$ has the following dual graph 
$$(-2) -- (-1) -- (-2) -- (-2).$$ 
Here $D_1$ is the proper transform of $L_{P_1P_2}$. 
Let ${\widetilde Y} \rightarrow Y$ be the contraction of $D = \sum_{i=1}^3 D_i$. 
Then this $Y$ is the {\it unique} Gorenstein del Pezzo surface 
of Picard number 1 and with $\text{\rm Sing} Y = A_2 + A_1$. Note that $E, D_i$ are 
the only negative curves on ${\widetilde Y}$ [Y1, 2]. 
Clearly, $G := \text{\rm Aut}(Y) = \text{\rm Aut}({\widetilde Y})$. Since $G$ clearly stabilizes 
all negative curves, there is an induced $G$-action on ${\bold P}^2$ such that 
the map $\mu$ is $G$-equivariant. 
Thus $G = \{g \in PGL_2({\bold C})| g(P_1) = P_1, g(L_{P_1P_2}) = L_{P_1P_2}\}$. 
So $G = \{(a_{ij}) \in PGL_2({\bold C}) | a_{ij} = 0 (i < j)\}$. 
Since $D_1 \cong {\bold P}^1$, each element of $G$ of finite order 
fixes either two or all points of $D_1$. 

\par \vskip 1pc \noindent 
{\bf Claim 1.6.1.} {\it Let $p \ge 2$ be a prime number. 
Given any point $Q \in D_1 \setminus \{D_1 \cap E\}$, 
there are exactly $p+1$ different ${\bold Z}/(p)$-actions on ${\widetilde Y}$, 
modulo equivariant isomorphism, 
each of which fixes at least two points of $D_1$: $Q$, $D_1 \cap E$; 
exactly one of these $p+1$ actions fixes three and hence all points of 
$D_1$.}

\par \vskip 1pc 
We use the same $Q$ to denote its image on ${\bold P}^2$ 
and we may assume that $Q = P_2$ after change of coordinates. 
Now $g = (a_{ij})$ in $G$ fixes $P_2$ if and only if $a_{32} = a_{ij} = 0$ ($i < j$). 
We may assume that $a_{11} = 1$ and write $g = g_1 + h$ 
where $g_1 = \text{\rm diag}[1, a_{22}, a_{33}]$ and $h = (h_{ij})$ 
with $h_{i1} = a_{i1}$ for for $i = 2, 3$ and 
all other $h_{ij} = 0$. Then $h^2 = 0$. So if $g$ is of prime order $p$,
then $I_3 = g^p = (g_1 + h)^p = g_1^p + p g_1^{p-1}h$. 
This equality implies that $g_1^p = I_3$ and $h = 0$. 
So $\text{\rm ord}(g) = p$ if and only if $\langle g \rangle =
\langle g_i \rangle$ for some $0 \le i \le p$, where  
$g_p = \text{\rm diag}[1, 1, \zeta_p]$ 
and $g_b = \text{\rm diag}[1, \zeta_p, \zeta_p^b]$ ($0 \le b \le p-1$)
with $\zeta_p = exp(2 \pi \sqrt{-1}/p)$. Now the image 
$\mu(D_1) = L_{P_1P_2} = \{X = 0\}$ is $g$-fixed if and only if 
$g = [1, \zeta_p, \zeta_p]$. The claim follows. 

\par
For an arbitray element $g$ of prime order $p$ in $G$,
we know that $g$ fixes at least two distinct points $D_1 \cap E$ and $Q$ on $D_1$.
Clearly, there is a projective
tranformation $\tau$ mapping $P_1, P_2$ to $P_1$, (the image of) $Q$.
This $\tau$ lifts to an automorphism $\tau$ on ${\widetilde Y}$ mapping
(the pre-image of) $P_2$ to $Q$.
Now $\langle \tau^{-1}g\tau \rangle$ fixes $P_1$ and $P_2$
and hence equals one of the $\langle g_i \rangle$ above for some $i$.
So modulo equivariant isomorphism, $\langle g_i \rangle$ ($0 \le i \le p$)
are the only non-trivial actions of ${\bold Z}/(p)$ on $Y$.

\par \vskip 1pc \noindent 
{\bf Example 1.7.} Let $Y$ be the {\it unique} Gorenstein del Pezzo 
surface of Picard number 1 and with $\text{\rm Sing} Y = A_3+2A_1$. Then $K_Y^2 = 4$. 

\par 
We first construct such unique surface $Y$. 
Let $P_1 = [0, 0, 1]$, $P_2 = [0, 1, 0]$ and $P_3 = [1, 0, 0]$. 
Let $\mu: {\widetilde Y} \rightarrow {\bold P}^2$ 
be the blow up of $P_1, P_2$ and 3 infinitely near 
points of them such that 
$\mu^{-1}(L_{P_1P_2} + L_{P_1P_3}) = D_1 + E_1 + D_4 + D_3 + D_2 + E_2 + D_5$ 
has the following dual graph 
$$(-2) -- (-1) -- (-2) -- (-2) -- (-2) -- (-1) -- (-2).$$ 
Here $D_1, D_2$ are the proper transforms of $L_{P_1P_3}$, $L_{P_1P_2}$. 
Let ${\widetilde Y} \rightarrow Y$ be the contraction of $D = \sum_{i=1}^5 D_i$. 
Then this $Y$ is the {\it unique} Gorenstein del Pezzo surface of Picard number 1  
and with $\text{\rm Sing} Y = A_3 + 2A_1$. Note that $E_i, D_j$ 
are the only negative curves on ${\widetilde Y}$ 
[Y1, 2]. Clearly, $G := \text{\rm Aut}(Y) = \text{\rm Aut}({\widetilde Y})$. 
Let $H = \{h \in G | h(D_1) = D_1\}$. Then $H$ stabilizes every negative curve. 
As in Example 1.6, $\mu$ is $H$-equivariant and 
$H = \{h \in PGL_2({\bold C})| h(P_i) = P_i (i=1,2), h(L_{P_1P_3}) = L_{P_1P_3}\}$. 
So $H = \{(a_{ij}) \in PGL_2({\bold C}) | a_{ij} = 0$ if $i \ne j$ 
and $(i, j) \ne (3,1)\}$. 

\par \vskip 1pc \noindent 
{\bf Claim 1.7.1.} $G \cong H \rtimes {\bold Z}/(2)$. 

\par \vskip 1pc 
Indeed, if $\sigma$ is in $G$ but not in $H$, then $\sigma$ switchs $D_1$ and $D_5$ 
and $\sigma^2$ is in $H$. On the other hand, the blow down 
${\widetilde Y} \rightarrow {\bold P}^1 \times {\bold P}^1$ 
of $E_1+D_4$ and $E_2+D_2$ to points $(0, 0), (\infty, 0)$ is $G$-equivariant 
and the automorphism $(x, y) \mapsto (1/x, y)$ downstairs lifts 
to an involution $\sigma$ in $G$ switching $D_1$ and $D_5$. 

\par \vskip 1pc 
By the claim above, every prime order $p$ ($p \ge 3$) element $g$ of $G$ 
is contained in $H$. As in Example 1.6, 
$\langle g \rangle = \langle \text{\rm diag}[1, 1, \zeta_p] \rangle$ 
or $\langle g \rangle = \langle \text{\rm diag}[1, \zeta_p, \zeta_p^b] \rangle$ ($0 \le b \le p-1$).
 
\par \vskip 1pc \noindent 
{\bf Example 1.8.} Let $Y$ be the {\it unique} Gorenstein del Pezzo 
surface of Picard number 1 and with $\text{\rm Sing} Y = 3A_2$. Then $K_Y^2 = 3$. 

\par 
We first construct such unique surface $Y$. 
Let $P_1 = [0, 0, 1]$, $P_2 = [0, 1, 0]$ and $P_3 = [1, 0, 0]$. 
Let $\mu: {\widetilde Y} \rightarrow {\bold P}^2$ 
be the blow up of $P_1, P_2, P_3$ and 3 infinitely near 
points of them such that 
$\mu^{-1} (\sum_{i < j} L_{P_iP_j}) = 
E_1 + D_1 + D_2 +  E_2 + D_3 + D_4 + E_3 + D_5 + D_6$ 
is a simple loop (with $E_1 . D_6 = 1$ and) with 
$E_i^2 = -1$ and $D_j^2 = -2$. 
Here $D_1, D_3, D_5$ are the proper transforms of $L_{P_iP_j}$ 
with $(i, j) = (1, 2)$, $(1, 3)$, $(2, 3)$. 
Let ${\widetilde Y} \rightarrow Y$ be the contraction of $D = \sum_{i=1}^6 D_i$. 
Then this $Y$ is the {\it unique} Gorenstein del Pezzo surface of Picard number 1  
and with $\text{\rm Sing} Y = 3A_2$. Note that $E_i, D_j$ 
are the only negative curves on ${\widetilde Y}$ 
[Y1, 2]. Clearly, $G := \text{\rm Aut}(Y) = \text{\rm Aut}({\widetilde Y})$. 
Let $H = \{h \in G | h(E_i) = E_i (i=1,2,3)\}$. 
Then $H$ stabilizes every negative curve. 
As in Example 1.6, the $\mu$ is $H$-equivariant and 
$H = \{h \in PGL_2({\bold C})| h(P_i) = P_i (i=1,2, 3)\}$. 
So $H = \{\text{\rm diag}[a, b, c] \in PGL_2({\bold C})\}$. 

\par \vskip 1pc \noindent 
{\bf Claim 1.8.1.} $G \cong H \rtimes S_3$. 

\par \vskip 1pc 
Indeed, if $\sigma$ is in $G$ but not in $H$, then $\sigma$ 
permutes $E_i$'s and hence $G/H \le S_3$. 
On the other hand, consider the blow down 
${\widetilde Y} \rightarrow {\bold P}^1 \times {\bold P}^1$ 
of $E_1+D_6$, $E_3$ and $E_2+D_3$ to points 
$(0, 0)$, $(\infty, 0)$ and $(\infty, \infty)$. 
Then the involution $\sigma : (x, y) \mapsto (1/y, 1/x)$ 
downstairs lifts to an involution $\sigma$ in $G$ switching $E_1$ and $E_2$. 
Similarly, we can find an involution $\sigma_{ij}$ in $G$ 
switching $E_i$ and $E_j$. This proves the claim. 

\par \vskip 1pc 
By the claim above, every prime order $p$ ($p \ge 5$) element $g$ of $G$ 
is contained in $H$. As in Example 1.6, 
$\langle g \rangle = \langle \text{\rm diag}[1, 1, \zeta_p] \rangle$ 
or $\langle g \rangle = \langle \text{\rm diag}[1, \zeta_p, \zeta_p^b] \rangle$ ($0 \le b \le p-1$). 

\par \vskip 1pc \noindent 
{\bf Example 1.9.} Let $Y$ be the {\it unique} Gorenstein del Pezzo 
surface of Picard number 1 and with $\text{\rm Sing} Y = 2A_3+A_1$. Then $K_Y^2 = 2$. 

\par 
We first construct such unique surface $Y$. 
Let $P_1 = [0, 0, 1]$, $P_2 = [0, 1, 0]$, $P_3 = [1, 0, 0]$ 
and $P_4 = [1, 1, 0]$. 
Let $\mu: {\widetilde Y} \rightarrow {\bold P}^2$ 
be the blow up of $P_1, P_2, P_3$ and 4 infinitely near 
points of them such that 
$\mu^{-1} (\sum_{i < j} L_{P_iP_j}) = L + C$ with 
$L = E_1 + D_1 + D_2 + D_3 + E_4 + D_6 + D_5 + D_4$ and $C = E_2 + D_7 + E_3$. 
Here $E_i^2 = -1$ and $D_j^2 = -2$. 
This $L$ is a simple loop (with $E_1 . D_4 = 1$) 
and $C$ is a linear chain such that $L . C = 2$ and $L$ 
meets $C$ at the two points $E_2 \cap D_2$ and $E_3 \cap D_5$. 
Also $D_1, D_5, D_6, E_2$ are the proper transforms of $L_{P_iP_j}$ 
with $(i, j) = (1, 2)$, $(2, 3)$, $(1, 4)$, $(1, 3)$. 
Let ${\widetilde Y} \rightarrow Y$ be the contraction of $D = \sum_{i=1}^7 D_i$. 
Then this $Y$ is the {\it unique} Gorenstein del Pezzo surface of Picard number 1
and with $\text{\rm Sing} Y = 2A_3+A_1$. Note that $E_i, D_j$ 
are the only negative curves on ${\widetilde Y}$ 
[Y1, 2]. Clearly, $G := \text{\rm Aut}(Y) = \text{\rm Aut}({\widetilde Y})$. 
Let $H = \{h \in G | h(E_i) = E_i (1 \le i \le 4)\}$. 
Then $H$ stabilizes every negative curve. 
As in Example 1.6, $\mu$ is $H$-equivariant and 
$H = \{h \in PGL_2({\bold C})| h(P_i) = P_i (1 \le i \le 4)\}$. 
So $H = \{\text{\rm diag}[1, 1, c] \in PGL_2({\bold C})\}$. 

\par \vskip 1pc \noindent 
{\bf Claim 1.9.1.} $G \cong H \rtimes ({\bold Z}/(2))^{\oplus 2}$. 

\par \vskip 1pc
Indeed, if $\sigma$ is in $G$ but not in $H$, then $\sigma$ 
stabilizes both sets $\{E_1, E_4\}$ and $\{E_2, E_3\}$. 
Note that the blow down 
${\widetilde Y} \rightarrow {\bold P}^1 \times {\bold P}^1$ 
of $E_1$, $E_4$, $E_2+D_2$ and $E_3+D_5$ to points 
$(0, \infty), (\infty, 0), (0, 0), (\infty, \infty)$ is $G$-equivariant. 
The involutions $\sigma_1 : (x, y) \mapsto (1/y, 1/x)$, 
$\sigma_2 : (x, y) \mapsto (y, x)$ 
lift to involutions $\sigma_1$, $\sigma_2$ in $G$ 
switching respectively $E_2$ and $E_3$, 
$E_1$ and $E_4$. This proves the claim. 

\par \vskip 1pc 
By the claim above, every prime order $p$ ($p \ge 3$) element $g$ of $G$ 
is contained in $H$. Arguing as in Example 1.6, 
$\langle g \rangle = \langle \text{\rm diag}[1, 1, \zeta_p] \rangle$. 

\par \vskip 1pc \noindent 
{\bf Example 1.10.} Let $Y$ be the {\it unique} Gorenstein del Pezzo 
surface of Picard number 1 and with $\text{\rm Sing} Y = A_7$. Then $K_Y^2 = 2$.
We shall show that $\text{\rm Aut}(Y) = {\bold Z}/(4) \times {\bold Z}/(2)$. 

\par 
We first construct such unique surface $Y$. 
Let $P_1 = (0, 0)$ and $P_2 = (\infty, 0)$.
Let $\mu: {\widetilde Y} \rightarrow {\bold P}^1 \times {\bold P}^1$ 
be the blow up of $P_1, P_2$ and 4 infinitely near 
points of them such that 
$\mu^{-1} (L_{P_1} + M + L_{P_2}) = E_1 + D + E_2$ with 
$D = D_1 + \cdots + D_7$; here $L_{P_i}$ ($i = 1,2$) are the fibres
$x = 0$ and $x = \infty$, and $M$ is the section $y = 0$.
Also $E_i^2 = -1$ and $D_j^2 = -2$. 
This $D$ is a linear chain
such that $E_i . D = E_i . D_j = 1$ with $(i, j) = (1, 2)$, $(2, 6)$.
Also $D_1$, $D_4$, $D_7$ are the proper transforms of $L_{P_1}$, $M$, $L_{P_2}$.
Let ${\widetilde Y} \rightarrow Y$ be the contraction of $D$. 
Then this $Y$ is the {\it unique} Gorenstein del Pezzo surface of Picard number 1
and with $\text{\rm Sing} Y = A_7$. Note that $E_i, D_j$ 
are the only negative curves on ${\widetilde Y}$ 
[Y1, 2]. 

\par
Clearly, $G := \text{\rm Aut}(Y) = \text{\rm Aut}({\widetilde Y})$. 
Let $H = \{h \in G | h(E_i) = E_i (i = 1,2)\}$. 
Then $H$ stabilizes every negative curve. 
Note that the restriction $H|D_i$ with $i = 2$ (resp. $i = 6$)
is trivial for $H$ fixes three intersection points of $D_i$ ($\cong {\bold P}^1$)
with $D_{i-1}$, $D_{i+1}$ and $E_1$ (resp. $E_2$).
As in Example 1.6, $\mu$ is $G$-equivariant.
So $H = \{h \in \text{\rm Aut}({\bold P}^1 \times {\bold P}^1) | h(P_i) = P_i (i=1,2);
h|D_j = \text{\rm id} (j=2,6)\}$. 

\par
Since an element $h$ in $H$
stabilizes the curves $L_{P_1} = \{x = 0\}$
and $M = \{y = 0\}$, we have $h : (x, y) \mapsto (ax, by)$.
Following the blow up process,
we see that $h|D_j = \text{\rm id}$ ($j = 2, 6$) if and only if $a = b^2$ and $a^{-1} = b^2$.
Thus $H = \langle h_1 \rangle \cong {\bold Z}/(4)$,
where $h_1 = (-1, \sqrt{-1})$.
Let $\sigma_c : (x, y) \mapsto (c/x, y)$ be an involution.
Then for a unique choice of $c$, the $\sigma_c$ lifts
to an involution on ${\widetilde Y}$ switching $E_1$ and $E_2$.
Clearly, $h_1 \sigma_c = \sigma_c h_1$.
Therefore, $G = \langle h_1 \rangle \times \langle \sigma_c \rangle 
\cong {\bold Z}/(4) \times {\bold Z}/(2)$.

\head 
Section 2. Proofs of Theorems 
\endhead 

We will first prove Theorems 2.1-2.3. 
Theorems A, B and C will follow from Theorem 2.3,
the observation that $K_Y^2 = 9 - \#D$ in Proposition 2.5,
Proposition 2.5 (3) and Lemma 2.6 (4).

\par \vskip 1pc \noindent 
{\bf Theorem 2.1}. {\it Let $p \ge 5$ be a prime number. 
Let $f : X \rightarrow {\bold P}^1$ be an extremal rational elliptic (smooth) surface 
such that ${\bold Z}/(p)$ acts non-trivially on $X$. 
Then modulo equivariant isomorphism, the ${\bold Z}/(p)$-action equals 
one of those in Examples $1.1 - 1.5$ ($4$ isolated ones and a family parametrized 
by a parameter $J$ in ${\bold C}$).
In particular, ${\bold Z}/(q) \subseteq \text{\rm Aut}(X)$ for every
prime $q \ge 5$, unless $p = 5$ and $X$ is given in Example $1.5$.} 

\par \vskip 1pc \noindent 
{\bf Theorem 2.2.} {\it Let $p \ge 5$ be a prime number. 
Let $f : X \rightarrow {\bold P}^1$ be an extremal rational elliptic (smooth) surface. 
Then a non-trivial action (if exists) of ${\bold Z}/(p)$ on $X$ 
is unique modulo equivariant isomorphism.
Also the action stabilizes all negative curves on $X$, unless
$p = 5$ and $X$ is given in Example $1.5$.} 

\par \vskip 1pc \noindent 
{\bf Theorem 2.3.} {\it Let $p \ge 5$ be a prime number. 
Let $Y$ be a Gorenstein del Pezzo (singular) surface of Picard number $1$.} 

\par \noindent 
(1) {\it $Y$ admits a non-trivial ${\bold Z}/(p)$-action 
if and only if either $K_Y^2 \ge 2$ and $\text{\rm Sing} Y \ne A_7$, or $K_Y^2 = 1$
and $|-K_Y|$ has exactly two singular members
(i.e., $K_Y^2 = 1$ and $Y$ equals one of those four in
Examples $1.1 \sim 1.4$)
whence $\text{\rm Sing} Y$ is one of the following (see {\bf (2.5)}.)} 
$$E_8, \, E_7 + A_1, \, E_6 + A_2, \, 2D_4.$$ 

\par \noindent 
(2) {\it If $K_Y^2 = 1$, then there is at most one non-trivial 
${\bold Z}/(p)$-action 
on $Y$, modulo equivariant isomorphism.} 

\par \noindent 
(3) {\it If $\text{\rm Sing} Y$ is either $A_2+A_1$, or $A_3+2A_1$ or $3A_2$, 
then there are exactly $p+1$ different non-trivial ${\bold Z}/(p)$-actions 
on $Y$, modulo equivariant isomorphism (see Examples $1.6$, $1.7$ and $1.8$).} 

\par \noindent
(4) {\it If $\text{\rm Sing} Y = A_7$, then $\text{\rm Aut}(Y) \cong
{\bold Z}/(4) \oplus {\bold Z}/(2)$ (see Example $1.10$).}

\par \noindent 
(5) {\it If $K_Y^2 \ge 2$, $Y \ne \overline{\Sigma}_2$ 
and $\text{\rm Sing} Y$ is neither one of those four in $(3)$ and $(4)$, 
then there is a unique non-trivial ${\bold Z}/(p)$-action 
on $Y$, modulo equivariant isomorphism.} 

\par \vskip 1pc 
We need some preparations. 

\par \vskip 1pc \noindent 
{\bf Lemma 2.4.} {\it Let $f : X \rightarrow B (\cong {\bold P}^1)$
be an extremal rational elliptic (smooth) surface such that there are at least
three singular fibres.
Then either $|\text{\rm Aut}(X)| = 2^a 5$ with $1 \le a \le 3$ or
$|\text{\rm Aut}(X)| = 2^a3^b$ with $1 \le a+b \le 7$. 
In the first case every element of order $5$ acts transitively
on the set of the $5$ sections of $f$.}

\par \vskip 1pc \noindent
{\it Proof.} Let $F$ be the generic fibre of $f$ over the function field
${\bold C}(B)$. 
We have an exact sequence (see [G, p.128]):
$$(1) \rightarrow \text{\rm Aut}(F) \rightarrow \text{\rm Aut}(X) \rightarrow \text{\rm Aut}(B).$$
Let $F_1 = f^{-1}(t_1), \dots, F_r = f^{-1}(t_r)$ be 
all singular fibres of $f$.
Then $H = \text{\rm Im}(\text{\rm Aut}(X) \rightarrow \text{\rm Aut}(B))$
acts on the set $\{t_1, \dots, t_r\}$.
The natural map $H \rightarrow S_r$ is injective for $r \ge 3$
and $B \cong {\bold P}^1$.
If we divide $\{t_1, \dots, t_r\}$ into a disjoint union 
of subsets of cardinality $n_k$ with $\sum_k n_k = r$,
such that fibres over points in the same subset are of the same type,
then $H$ stabilizes each subset and hence $|H|$ divides $\Pi_k (n_k)!$.

\par
On the other hand, if we let $\text{\rm MW}(f)$ be the Mordell-Weil
group of all sections of $f$ and
$\text{\rm Aut}(F)_0$ be the subgroup of $\text{\rm Aut}(F)$
fixing the zero element (coming from a pre-designated zero section of $f$),
then $\text{\rm Aut}(F) = \text{\rm MW}(f) \rtimes \text{\rm Aut}(F)_0$. 
So we conclude that $|\text{\rm Aut}(X)| = |\text{\rm Aut}(F)_0||\text{\rm MW}(f)||H|$.
All such $f$ are given in [MP, Theorem 4.1] (see also [ibid, Theorem 5.4]).
By [ibid. Table 5.1-5.3], the $J$-function $J(F)$ is not a constant.
So $|\text{\rm Aut}(F)_0| = 2$. Now the lemma follows from the classification
of singular fibre types and $\text{\rm MW}(f)$ in [ibid. Theorem 4.1]
and the above observation that $|H|$ divides $\Pi_k (n_k)!$.

\par \vskip 1pc \noindent
{\bf Proposition 2.5.} 
{\it Let $Y$ ($\ne {\overline \Sigma}_2$) 
be a Gorenstein del Pezzo (singular) surface of Picard number $1$.} 

\par \noindent 
(1) {\it $\text{\rm Sing} Y$ is one of the following ($26$ of them):} 
$$\gather 
E_8, \, E_7+A_1, \, E_6+A_2, \, 2D_4, \, D_8, \, 
D_5+A_3, \, D_6+2A_1, \, A_8, \, A_7+A_1, \\ 
2A_4, \, A_5+A_2+A_1, \, 2A_3+2A_1, \, 4A_2, \, 
E_7, \, D_6+A_1, \, A_7, \, A_5+A_2, \, D_4+3A_1, \\ 
2A_3+A_1, \, 3A_2, \, A_5+A_1, \, E_6, \, 
A_3+2A_1, \, D_5, \, A_4, \, A_2+A_1. 
\endgather$$ 

\par \noindent 
(2) {\it $\text{\rm Sing} Y$ determines uniquely the isomorphism class of $Y$
and $|-K_Y|$ has at least three singular members, if 
$\text{\rm Sing} Y$ is not $E_8, E_7+A_1, E_6+A_2$ or $2D_4$. 
If $D$ is one of $E_8, E_7+A_1$ and $E_6+A_2$, then there are 
exactly two isomorphism classes of $Y = Y_D(i)$ with $\text{\rm Sing} Y = D$;
$Y = Y_D(1)$ is given in Examples $1.1 - 1.3$; the
$|-K_Y|$ with $Y = Y_D(1)$ (resp. $Y = Y_D(2)$) has exactly two 
(resp. at least three) singular members.
If $\text{\rm Sing} Y = 2D_4$ then $Y$ is isomorphic to one of those 
$Y_{J, 0}$ in Example $1.4$ parametrized by $J$ in ${\bold C}$
and $|-K_Y|$ has exactly two singular members.}
\par \noindent
(3) {\it Suppose that $K_Y^2 = 1$ and $|-K_Y|$ contains at least 
three singular members. Then $|\text{\rm Aut}(Y)| = 2^a3^b$ with $1 \le a+b \le 7$.}

\par \vskip 1pc
For the $Y$ in the Proposition above, we let 
Let ${\widetilde Y} \rightarrow Y$ 
be the minimal resolution and $D$ the exceptional divisor.
Clearly, $\text{\rm Aut}(Y) = \text{\rm Aut}({\widetilde Y})$. 
Also the Picard number
$\rho({\widetilde Y}) = 1 + \#D$ and 
$K_Y^2 = K_{\widetilde Y}^2 = 9 - \#D$, where 
$\#D$ is the number of irreducible components of $D$. 

\par \vskip 1pc \noindent
{\it Proof of Propositon $2.5$.} (1) is well known (see e.g. [MZ1]).
(2) is proved in [Y1,2]. We state the idea of the proof of (2)
for the case $K_Y^2 = 1$. 
Let $X \rightarrow {\widetilde Y}$ be the blow up of the unique point
$\text{\rm Bs}|-K_{\widetilde Y}|$ with $E$ the exceptional curve
(see the proof of Theorem 2.3 below).
Then $f = \Phi_{|-K_X|} : X \rightarrow B (\cong {\bold P}^1)$
is an extremal rational elliptic surface. These $X$, $\text{\rm MW}(f)$
and isomorphism classes are classified in [MP, Theorems 4.1 and 5.4].
The composition map $X \mapsto Y$
defines a one-to-one correspondence between the set of
isomorphism classes of extremal rational elliptic (smooth) surfaces and
the set of Gorenstein del Pezzo surfaces of Picard number 1 and degree 1;
the singular members of $|-K_Y|$ are exactly
the images of singular fibres of $f$.

\par
The assertion(3) follows from Lemma 2.4 and the observation
that $\text{\rm Aut}(Y) = \{g \in \text{\rm Aut}(X) | g(E) = E)\} \supseteq \langle \sigma \rangle$; here
$E$ is a section of $f$ and $\sigma$ acts as $-\text{\rm id}$ on the generic fibre
(fixing the zero section which is chosen to be $E$).

\par \vskip 1pc \noindent 
{\bf Lemma 2.6.} {\it Let $Y$ be as in the Proposition above. 
Let $g$ be an element in $\text{\rm Aut}(Y) = \text{\rm Aut}({\widetilde Y})$
of prime order $p$ ($p \ge 5$).} 

\par \noindent 
(1) {\it Every negative curve on ${\widetilde Y}$ is 
either a $(-1)$-curve or a $(-2)$-curve, and all negative curves
are $g$-stable.
If $K_Y^2 \ge 2$, then there are at most four of $(-1)$-curves 
on ${\widetilde Y}$.}

\par \noindent
(2) {\it If $E$ is a $(-1)$-curve on ${\widetilde Y}$, then 
$-K_{\widetilde Y} = dE + \Delta$ numerically, where $d = K_Y^2$ and 
$\Delta$ is an effective ${\bold Q}$-divisor with support in $D$.} 

\par \noindent
(3) {\it There is a $(-1)$-curve $E$ on ${\widetilde Y}$
such that $E . D = 1$ if $\text{\rm Sing} Y$ equals one of the following:}
$$E_7, \, D_6+A_1, \, A_7, \, A_5+A_2, \, A_5+A_1, \, E_6, \, D_5, \, A_4.$$

\par \noindent
(4) {\it Theorem A $(3)$ in the Introduction is true.}

\par \vskip 1pc \noindent
{\it Proof.} The nef and bigness of $-K_{\widetilde Y}$ implies that
every negative curve is either a $(-1)$-curve or a $(-2)$-curve.
When $K_Y^2 = 1$, (1) is reduced to Theorem 2.2 
(see the proof of Theorem 2.3). 
(3) and the second part of (1) follows from Figure 5 in [Y1] or [Y2, Ch 4].
In particular, when $K_Y^2 \ge 2$ all $(-1)$-curve are $g$-stable. Also every
component of $D$ of type as shown in Proposition 2.5 is clearly $g$-stable.
So (1) follows. 

\par
For (2), we note that the Picard number $\rho(Y) = 1$,
whence we have numerically $-K_Y = b \overline{E}$ for some rational number
$b$, where $\overline{E}$ is the image on $Y$ of $E$.
Pulling back this relation,
we get $-K_{\widetilde Y} = b(E + \Delta')$, where $\Delta'$ is 
an effective ${\bold Q}$-divisor
and supported by $D$. Intersecting this with $-K_{\widetilde Y}$,
we get $b = K_{\widetilde Y}^2 = K_Y^2$. 

\par
To prove (4), we need:

\par \vskip 1pc \noindent
{\bf Claim 2.6.1.} {\it There is a birational morphism 
$\mu : {\widetilde Y} \rightarrow {\bold P}^2$.}

\par \vskip 1pc
By (1), we have a birational morphism ${\widetilde Y} \rightarrow \Sigma_d$
where $\Sigma_d$ is the Hirzebruch surface of degree $d \le 2$.
We let $M$ be the section on $\Sigma_d$ with $M^2 = -d$.
Since $K_{\widetilde Y}^2 < 8$, the map
${\widetilde Y} \rightarrow \Sigma_d$ factors through the
blow up $Z \rightarrow \Sigma_d$ of a point on a fibre $F$.
The inverse on $Z$ of $F$ is a pair $E+F'$ of intersecting $(-1)$-curves
with $F'$ the proper transform of $F$. 
If $d = 1$, then the claim is clear. 
If $d = 0, 2$, we have a composition
${\widetilde Y} \rightarrow Z \rightarrow {\Sigma}_1$
where the latter map is the blow down of $F'$.
The claim follows.

\par \vskip 1pc
We continue the proof of Lemma 2.6 (4).
By (1), the $\mu$-exceptional divisor is $g$-stable.
So there is an induced action of $g$ on ${\bold P}^2$
so that $\mu$ is $g$-equivariant. Now Theorem A (3) follows.
This proves the lemma.

\par \vskip 1pc
Let $X$ be a smooth projective surface with an 
order-$n$ automorphism $g$. Let $C_0 + C_1 + \cdots + C_r$ be a linear 
chain of $(-2)$-curves each of which is $g$-stable. 
Denote by $P_{i+1} = C_ i \cap C_{i+1}$. Note that 
the fixed set $C_i^g$ contains two distinct points $P_i$, $P_{i+1}$, 
where $P_0 \in C_0$ and $P_{r+1} \in C_r$. Let $\zeta$ be a primitive $n$-th 
root of 1. 

\par \vskip 1pc \noindent 
{\bf Lemma 2.7.} {\it Let $g|P_i = (\zeta^{a_i}, \zeta^{b_i})$ be the diagonalization. 
Then we have $b_i + a_{i+1} = 0$ (\text{\rm mod} $n$) and $a_i + b_i = a_j + b_j$ (\text{\rm mod} $n$) 
For all $i$ and $j$.} 

\par \vskip 1pc \noindent 
{\it Proof.} The first equality follows from the fact if $g$ acts as a multiple $\zeta^e$ 
at the origin of $C_i \cong {\bold P}^1$ then it acts as a multiple $\zeta^{-e}$ 
at the infinity. For the second equality, 
we note that the restricted line bundle ${\Cal O}(K_X)|C_i$ is trivial 
for $C_i$ is a $(-2)$-curve, and $g$ acts on it by a constant 
multiple $g : (dx_i \wedge dy_i) |C_i  \mapsto \zeta^{a_i+b_i} 
(dx_i \wedge dy_i) |C_i$, where $x_i, y_i$ are defining equations 
of $C_i, C_{i+1}$ at $P_i$. 
So $\zeta^{a_i+b_i} = 
\zeta^{a_{i+1} + b_{i+1}}$ for $P_i$ and $P_{i+1}$ are on the same curve $C_i$. 
This proves the lemma.

\par \vskip 1pc 
In Lemmas 2.8 and 2.9 below, we let 
$f : X \rightarrow B \cong {\bold P}^1$ be an extremal 
rational elliptic (smooth) surface and $g$
an element in $\text{\rm Aut}(X)$ of prime order $p$ ($p \ge 5$).
Note that each fibre is a member of $|-K_X|$. So the fibration is $g$-stable. 
We denote by $g|B$ the induced action of $g$ on the base curve $B$. 

\par \vskip 1pc \noindent 
{\bf Lemma 2.8.} {\it Suppose that $g|B = \text{\rm id}$. Then $f$ has singular fibres of types 
$I_5, I_5, I_1, I_1$ (such $X$ is unique; see $[MP]$), $p = 5$ 
and $\text{\rm MW}(f) \cong {\bold Z}/(5)$. 
Moreover, the action of $\langle g \rangle \cong {\bold Z}/(5)$ 
on $X$ is identical to that of $\text{\rm MW}(f)$ on $X$ as translations 
(on general fibre).} 

\par \vskip 1pc \noindent 
{\it Proof.} If $p$ does not divide $|\text{\rm MW}(f)|$, then $g$ stabilizes each section; 
so $g$ fixes each section for $g|B = \text{\rm id}$. Thus a general fibre $F$ is an elliptic curve 
on which $g$ acts with a fixed point (the intersection of $F$ and a section). 
So $\text{\rm ord}(g|F)$ divides 12. On the other hand $\text{\rm ord}(g|F)$ divides $\text{\rm ord}(g|X) = p \ge 5$. 
So $g|F = \text{\rm id}$. This and $g|B = \text{\rm id}$ imply that $g = \text{\rm id}_X$, a contradiction. 
Thus $p$ divides $|\text{\rm MW}(f)|$. Now the first part of lemma follows from 
[MP, Theorems 4.1 and 5.4]. Also $g$ acts transitively on the set 
of all 5 sections by the argument above. 

\par 
For the second part, choose an automorphism $h$ of $X$ 
coming from an element (denoted also by $h$) of $\text{\rm MW}(f)$ such that $g^{-1}h$ 
acts identically on the set of the 5 sections in $\text{\rm MW}(f)$. 
Then $k = g^{-1} h$ acts on a general fibre $F$ with at least 5 fixed points 
(the intersection of $F$ with the 5 sections). 
Note that $\text{\rm ord}(k|F)$ divides $12$ and if 
$k|F \ne \pm \text{\rm id}$, then $J(F) = 0$ or $1728$. 
On the other hand, $J(F)$ is not a constant by [MP, Table 5.3]. 
We reach a contradiction. So $k|F = \pm \text{\rm id}$.
If $k|F = \text{\rm id}$, then $k = \text{\rm id}_X$ for $k|B = \text{\rm id}$,
whence $g = h$. If $k|F = - \text{\rm id}$
then $g|F = (- \text{\rm id}_F) h$ has order $10$, a contradiction.
This proves the lemma. 

\par \vskip 1pc \noindent 
{\bf Lemma 2.9.} {\it Suppose that $g|B \ne \text{\rm id}$. Then $f$ has one of 
the following $4$ singular fibre types:} 

\par 
(i) $II, II^*, \,\,\,$ (ii) $III, III^*, \,\,\,$ 
(iii) $IV, IV^*, \,\,\,$ (iv) $I_0^*, I_0^*$. 

\par 
{\it Moreover, for each $k = 1, 2, 3$, there is a unique $X$ of singular fibre type $(k)$ 
above. The surfaces with singular fibre type $(iv)$ above are parametrized 
by the $J$-invariant of a general fibre $F$ ($J(F)$ is a constant).} 

\par \vskip 1pc \noindent 
{\it Proof.} If one singular fibre is not $g$-stable, then $f$ has $p$ 
($p \ge 5$) copies of 
the same fibre, which contradicts the fact that $f$ has at most four singular fibres 
[MP, Theorem 4.1]. So every singular fibre is $g$-stable. 
If $f$ has more than two singular fibres, then $g|B$ fixes more than two points 
(over which lying the singular fibres), so $g|B = \text{\rm id}$ for $B \cong {\bold P}^1$. 
This is a contradiction. Thus $f$ has at most two singular fibres. Now the lemma 
follows from [MP, Theorem 4.1]. 

\par \vskip 1pc 
In view of Lemmas 2.8 and 2.9, to show Theorems 2.1 and 2.2 we have only 
to show that for each $X$ satisfying one of 4 cases in Lemma 2.9, 
there is a unique non-trivial 
$G = {\bold Z}/(p)$-action ($p \ge 5$) on $X$ (one is given in Examples 1.1-1.4). 
Note that in all these 4 cases, $|\text{\rm MW}(f)| \le 4$ 
and hence $G$ acts trivially on the set of sections, i.e., 
$G$ stabilizes every section, because $p \ge 5$. This and $p \ge 5$ 
again imply that every component in the two singular fibres 
(each of which is $G$-stable) is $G$-stable 
and hence the central component $C_0$, meeting three other components $C_1, 
C_2, C_3$, of a singular fibre is $G$-fixed (point wise). 
Let $P = C_0 \cap C_1$ and let $h$ be a generator of $G$ such that $h|P$ can 
be diagonalized as $h|P = (1, \zeta_p)$, where $h$ acts identically 
along the direction of $C_0$ and as a multiple $\zeta_p = exp(2 \pi \sqrt{-1}/p)$ 
along the direction of $C_1$. 
Applying Lemma 2.7, we can show that the $g$ 
in Examples 1.1 - 1.4 satisfies $g|P = (1, \zeta_p)$ 
so that $g^{-1} h$ acts identically along the directions of $C_0$ and $C_1$. 
Hence $g^{-1} h = \text{\rm id}_X$ and $h = g$. This shows that the non-trivial action of 
${\bold Z}/(p)$ on $X$ is the same as that of $\langle g \rangle \cong {\bold Z}/(p)$ 
on $X$ in one of Examples 1.1 - 1.4. This proves Theorems 2.1 and 2.2.  

\par 
We can avoid the use of Lemma 2.7 by 
choosing a suitable power $g^e$ so that $g^e|P = (1, \zeta_p)$, whence $h = g^e$; 
so Lemma 2.7 claims that $e = 1$ indeed. 

\par \vskip 1pc 
Next we prove Theorem 2.3 by reducing to Theorems 2.1 and 2.2. 
Let ${\widetilde Y} \rightarrow Y$ and $D$
be as in Proposition 2.5. The assertion(4) is proved in Example 1.10.
If $K_Y^2 \ge 7$ then $K_Y^2 = 8$ and $Y$ equals the
quadric cone ${\overline \Sigma}_2$ in ${\bold P}^3$ (see [MZ1], or Proposition 2.5);
in this case Theorem 2.3 can be checked easily in the spirit of Example 1.6. 
So we may assume that $\text{\rm Sing} Y \ne A_7$ and $K_Y^2 \le 6$.

\par 
We first consider the case $K_Y^2 = 1$. 
By the Riemann-Roch theorem and Kawamata-Viehweg 
Vanishing [Ka, V], the linear system $|-K_{\widetilde Y}|$ has dimension 1, 
and it has a unique base point $Q$ by [D, Proposition 2, p.40]. 
Let $X \rightarrow {\widetilde Y}$ be the 
blow-up of $Q$ with $E$ the exceptional curve. Then the linear system 
$|-K_X|$ has dimension 1 and is base point free, so it defines an elliptic 
fibration $f : X \rightarrow B \cong {\bold P}^1$ with $E$ as a section. 
Clearly, the non-trivial action of ${\bold Z}/(p)$ on $Y$ induces 
a non-trivial action of ${\bold Z}/(p)$ 
on $X$ so that the birational morphism $X \rightarrow Y$ is ${\bold Z}/(p)$-equivariant. 
Since the Picard number $\rho(Y) = 1$, we see easily that $f$ is extremal 
in the sense of [MP]. Now Theorem 2.3 : the case $K_Y^2 = 1$
follows from Theorems 2.1 and 2.2. Indeed, since ${\bold Z}/(p)$ stabilizes the section $E$, 
this $X$ is not equal to the surface in Example 1.5. 

\par 
Next we consider the case $d := K_Y^2 \ge 2$. 
Suppose that $\text{\rm Sing} Y$ is as in Lemma 2.6 (3),
then there is a $(-1)$-curve $E_d$ on ${\widetilde Y}$ 
such that $E_d . D = 1$. If there is 
an element $g$ in $\text{\rm Aut}(Y) = \text{\rm Aut}({\widetilde Y})$
of prime order $p$ ($p \ge 5$),
then $g$ fixes at least two 
points on $E_d$ ($\cong {\bold P}^1$): $E_d \cap D$, $Q_d$. 
Let $\sigma : {\widetilde Y}_{d-1} \rightarrow {\widetilde Y}$ be the blow-up of $Q_d$
with $E_{d-1}$ the exceptional curve. 
Then the anti-canonical divisor of ${\widetilde Y}_{d-1}$ equals
$\sigma^*(-K_{\widetilde Y}) - E_{d-1} = dE_d' + (d-1)E_{d-1}  + \sigma^*(\Delta) \ge 0$
in notation of Lemma 2.6 (where $E_d'$ is the proper transform of $E_d$),
and hence it is nef and also big for its self-intersection equals $d-1 \ge 1$. 
So the divisor $D_{d-1} = \sigma^*(D) + E_d'$
of $(-2)$-curves which has zero intersection with the nef and big anti-canonical
divisor of ${\widetilde Y}_{d-1}$,
is contractible to rational double singularities by the map ${\widetilde Y}_{d-1}
\rightarrow Y_{d-1}$. Also $Y_{d-1}$ has Picard number 1 and ample $-K_{Y_{d-1}}$. 
Clearly, our $g$ in $\text{\rm Aut}({\widetilde Y}) = \text{\rm Aut}(Y)$ induces an element
$g$ in $\text{\rm Aut}({\widetilde Y}_{d-1}) = \text{\rm Aut}(Y_{d-1})$ so that $\sigma$ is $g$-equivariant.

\par
In particular, if $d = 2$ then $\text{\rm Sing} Y_{d-1}$ is one of those four in Theorem 2.3 (1)
and also such $\langle g \rangle \cong {\bold Z}/(p)$ in $\text{\rm Aut}(Y_{d-1})$ is unique modulo equivariant 
isomorphism; hence such $\langle g \rangle$ in $\text{\rm Aut}(Y)$ is also 
unique for any given $Y$ (see Lemma 2.6 (1)),
modulo equivariant isomorphism.
\par
For general $d \ge 2$, we have blow-ups below to reduce to the case $d = 1$:
${\widetilde Y}_1 \rightarrow {\widetilde Y}_2 \rightarrow 
\cdots \rightarrow {\widetilde Y}_d$.
Using the information on the location of the point $E_d \cap D$ in $D$
given in Figure 5 in [Y1] or [Y2, Ch 4],
all possible chains of singularity types :
$\text{\rm Sing} Y_1 \rightarrow \text{\rm Sing} Y_2 \rightarrow \cdots \rightarrow \text{\rm Sing} Y_d$
are as follows:
$$\gather
E_8 \rightarrow E_7 \rightarrow E_6 \rightarrow D_5 \rightarrow A_4, \\
E_7+A_1 \rightarrow D_6+A_1 \rightarrow A_5+A_1, \\
E_6+A_2 \rightarrow A_5+A_2.
\endgather $$
So Theorem 2.3 is true if $\text{\rm Sing} Y$ is one of those in Lemma 2.6 (3).

\par
We are left with those $Y$ with $\text{\rm Sing} Y$ equal to one of the following:
$$D_4+3A_1, \, 2A_3+A_1, \, 3A_2, \, A_3+2A_1, \, A_2+A_1.$$
For the last four cases, Theorem 2.3 follows from the arguments
in Examples 1.6-1.9.

\par
Now we consider the remaining case $\text{\rm Sing} Y = D_4+3A_1$.
Then ${\widetilde Y}$ contains a $(-1)$-curve $E_2$ meeting the three isolated
$(-2)$-curves in $D$ [ibid. Figure 5]. Since $g$ with $\text{\rm ord}(g) = p \ge 5$
stabilizes every component of $D$, this $g$ fixes the three points $D \cap E_2$
and hence $E_2$ ($\cong {\bold P}^1$) is $g$-fixed.
Blowing up a point on $E_2 \setminus D$, as above we will reduce to $Y_1$ with
$\text{\rm Sing} Y_1 = 2D_4$. Now the uniqueness of such 
$\langle g \rangle \cong {\bold Z}/(p)$ in $\text{\rm Aut}(Y)$
modulo equivariant isomorphism, follows from Theorem 2.3: the case $K_Y^2 = 1$.
This proves Theorem 2.3.

\head 
References 
\endhead 

\par \noindent
[dF] Tommaso de Fernex, Birational trasformations of prime order of the projective plane,
preprint 2001.

\par \noindent 
[D] M. Demazure, Lecture Notes in Mathematics, {\bf 777} (1980), Springer. 

\par \noindent
[DO] I. Dolgachev and D. Ortland, Point sets in projective spaces and theta functions,
Astérisque, Vol. {\bf 165} (1988).

\par \noindent
[G] M. H. Gizatullin, Rational $G$-surfaces,
Math. USSR Izv. {\bf 16} (1981), 103--134.

\par \noindent 
[GPZ] R. V. Gurjar, C. R. Pradeep, D. -Q. Zhang, 
On Gorenstein surfaces isomorphic to ${\bold P}^2/G$, Nagoya Math. J. to appear,
math.AG/0112242.

\par \noindent 
[H1] T. Hosoh, Automorphism groups of cubic surfaces, J. Algebra {\bf 192} (1997), 
651--677. 

\par \noindent 
[H2] T. Hosoh, Automorphism groups of quartic del Pezzo surfaces, 
J. Algebra {\bf 185} (1996), 374--389. 

\par \noindent
[I] V. A. Iskovskih, Minimal models of rational surfaces over 
arbitrary fields, Math. USSR Izv. {\bf 14} (1981), 17--39.

\par \noindent
[K] S. Kantor, Theorie der endlichen Gruppen von eindeutigen Transformationen
in der Ebene, Berlin: Mayer $\&$ Muller, 1895.

\par \noindent 
[Ka] Y. Kawamata, A generalization of Kodaira-Ramanujam's vanishing theorem, 
Math. Ann. {\bf 261} (1982), 43--46. 

\par \noindent 
[Ko] M. Koitabashi, Automorphism groups of generic rational surfaces, 
J. Algebra {\bf 116} (1988), 130 --142. 

\par \noindent
[KM] J. Kollar and S. Mori, Birational geometry of algebraic varieties,
Cambridge Tracts in Mathematics, {\bf 134} (1998).

\par \noindent
[M1] Yu. I. Manin, Rational surfaces over perfect fields, II, Math. USSR Sb.
{\bf 1} (1967), 141--168.

\par \noindent
[M2] Yu. I. Manin, Cubic forms : Algebra, geometry, arithmetic,
2nd ed, North-Holland Math. Library, {\bf 4} (1986),
North-Holland Publ. Co., Amsterdam-New York.

\par \noindent 
[MM] M. Miyanishi and K. Masuda, Open algebraic surfaces with finite group actions, 
Transform. Group, to appear. 

\par \noindent 
[MP] R. Miranda and U. Persson, On extremal rational elliptic surfaces, 
Math. Z. {\bf 193} (1986), 537--558. 

\par \noindent 
[MZ1,2] M. Miyanishi and D. -Q. Zhang, 
Gorenstein log del Pezzo surfaces of rank one, I; II, 
J. Algebra {\bf 118} (1988), 63--84; {\bf 156} (1993), 183--193. 

\par \noindent 
[MZ3] M. Miyanishi and D. -Q. Zhang, Equivariant classification of Gorenstein open 
log del Pezzo surfaces with finite group actions, Preprint 2001. 

\par \noindent
[O] K. Oguiso, Automorphism groups in a family of K3 surfaces,
math.AG/0104049.

\par \noindent 
[OS] K. Oguiso and T. Shioda, The Mordell-Weil lattice of a rational elliptic surface, 
Comment. Math. Univ. St. Paul. {\bf 40} (1991), 83--99. 

\par \noindent
[S] B. Segre, The Non-singular cubic surfaces, 
Oxford University Press, Oxford, 1942. 

\par \noindent 
[V] E. Viehweg, Vanishing theorems, J. Reine Angew. Math. {\bf 335} (1982), 1--8. 

\par \noindent 
[Y1] Q. Ye, On Gorenstein log del Pezzo surfaces, Japanese J. Math. to appear,
math.AG / 0109223.

\par \noindent 
[Y2] Q. Ye, On algebraic surfaces with non-positive Kodaira dimension, PhD thesis, 
National Univ. of Singapore, 2001. 

\par \noindent 
[ZD] D. -Q. Zhang, Automorphisms of finite order on rational surfaces, 
with an appendix by I. Dolgachev, J. Algebra {\bf 238} (2001), 560--589. 

\par \noindent 
[Z2] D. -Q. Zhang, Automorphisms of finite order
on extremal rational elliptic surfaces and 
Gorenstein del Pezzo surfaces of degree one,
Contemporary Math. to appear.

\enddocument